\documentclass[11pt]{article}
\usepackage{color}
\usepackage{booktabs}
\usepackage{booktabs,longtable}
\usepackage{amsfonts,amssymb,a4, graphicx, bm, makeidx}
\usepackage{amsmath}
\usepackage{soul}
\usepackage{mathrsfs}
\usepackage[all]{xy}
\usepackage[latin1]{inputenc}
\usepackage[dvips]{geometry}
\def\theequation{\thesection.\arabic{equation}}

\newcommand{\qed}{\hfill\rule{3mm}{3mm}}
\newtheorem{corollary}{Corollary}[section]
\newtheorem{theorem}[corollary]{Theorem}
\newtheorem{lemma}[corollary]{Lemma}

\newtheorem{proposition}[corollary]{Proposition}
\newtheorem{definition}[corollary]{Definition}

\makeatletter \@addtoreset{equation}{section} \makeatother
\definecolor{GreenYellow}{cmyk}{0.15,0,0.69,0}
\definecolor{Yellow}{cmyk}{0,0,1,0}
\definecolor{Goldenrod}{cmyk}{0,0.10,0.84,0}
\definecolor{Dandelion}{cmyk}{0,0.29,0.84,0}
\definecolor{Apricot}{cmyk}{0,0.32,0.52,0}
\definecolor{Peach}{cmyk}{0,0.50,0.70,0}
\definecolor{Melon}{cmyk}{0,0.46,0.50,0}
\definecolor{YellowOrange}{cmyk}{0,0.42,1,0}
\definecolor{Orange}{cmyk}{0,0.61,0.87,0}
\definecolor{BurntOrange}{cmyk}{0,0.51,1,0}
\definecolor{Bittersweet}{cmyk}{0,0.75,1,0.24}
\definecolor{RedOrange}{cmyk}{0,0.77,0.87,0}
\definecolor{Mahogany}{cmyk}{0,0.85,0.87,0.35}
\definecolor{Maroon}{cmyk}{0,0.87,0.68,0.32}
\definecolor{BrickRed}{cmyk}{0,0.89,0.94,0.28}
\definecolor{Red}{cmyk}{0,1,1,0}
\definecolor{OrangeRed}{cmyk}{0,1,0.50,0}
\definecolor{RubineRed}{cmyk}{0,1,0.13,0}
\definecolor{WildStrawberry}{cmyk}{0,0.96,0.39,0}
\definecolor{Salmon}{cmyk}{0,0.53,0.38,0}
\definecolor{CarnationPink}{cmyk}{0,0.63,0,0}
\definecolor{Magenta}{cmyk}{0,1,0,0}
\definecolor{VioletRed}{cmyk}{0,0.81,0,0}
\definecolor{Rhodamine}{cmyk}{0,0.82,0,0}
\definecolor{Mulberry}{cmyk}{0.34,0.90,0,0.02}
\definecolor{RedViolet}{cmyk}{0.07,0.90,0,0.34}
\definecolor{Fuchsia}{cmyk}{0.47,0.91,0,0.08}
\definecolor{Lavender}{cmyk}{0,0.48,0,0}
\definecolor{Thistle}{cmyk}{0.12,0.59,0,0}
\definecolor{Orchid}{cmyk}{0.32,0.64,0,0}
\definecolor{DarkOrchid}{cmyk}{0.40,0.80,0.20,0}
\definecolor{Purple}{cmyk}{0.45,0.86,0,0}
\definecolor{Plum}{cmyk}{0.50,1,0,0}
\definecolor{Violet}{cmyk}{0.79,0.88,0,0}
\definecolor{RoyalPurple}{cmyk}{0.75,0.90,0,0}
\definecolor{BlueViolet}{cmyk}{0.86,0.91,0,0.04}
\definecolor{Periwinkle}{cmyk}{0.57,0.55,0,0}
\definecolor{CadetBlue}{cmyk}{0.62,0.57,0.23,0}
\definecolor{CornflowerBlue}{cmyk}{0.65,0.13,0,0}
\definecolor{MidnightBlue}{cmyk}{0.98,0.13,0,0.43}
\definecolor{NavyBlue}{cmyk}{0.94,0.54,0,0}
\definecolor{RoyalBlue}{cmyk}{1,0.50,0,0}
\definecolor{Blue}{cmyk}{1,1,0,0}
\definecolor{Cerulean}{cmyk}{0.94,0.11,0,0}
\definecolor{Cyan}{cmyk}{1,0,0,0}
\definecolor{ProcessBlue}{cmyk}{0.96,0,0,0}
\definecolor{SkyBlue}{cmyk}{0.62,0,0.12,0}
\definecolor{Turquoise}{cmyk}{0.85,0,0.20,0}
\definecolor{TealBlue}{cmyk}{0.86,0,0.34,0.02}
\definecolor{Aquamarine}{cmyk}{0.82,0,0.30,0}
\definecolor{BlueGreen}{cmyk}{0.85,0,0.33,0}
\definecolor{Emerald}{cmyk}{1,0,0.50,0}
\definecolor{JungleGreen}{cmyk}{0.99,0,0.52,0}
\definecolor{SeaGreen}{cmyk}{0.69,0,0.50,0}
\definecolor{Green}{cmyk}{1,0,1,0}
\definecolor{ForestGreen}{cmyk}{0.91,0,0.88,0.12}
\definecolor{PineGreen}{cmyk}{0.92,0,0.59,0.25}
\definecolor{LimeGreen}{cmyk}{0.50,0,1,0}
\definecolor{YellowGreen}{cmyk}{0.44,0,0.74,0}
\definecolor{SpringGreen}{cmyk}{0.26,0,0.76,0}
\definecolor{OliveGreen}{cmyk}{0.64,0,0.95,0.40}
\definecolor{RawSienna}{cmyk}{0,0.72,1,0.45}
\definecolor{Sepia}{cmyk}{0,0.83,1,0.70}
\definecolor{Brown}{cmyk}{0,0.81,1,0.60}
\definecolor{Tan}{cmyk}{0.14,0.42,0.56,0}
\definecolor{Gray}{cmyk}{0,0,0,0.50}
\definecolor{Black}{cmyk}{0,0,0,1}
\definecolor{White}{cmyk}{0,0,0,0}

\begin{document}
\def\theequation{\thesection.\arabic{equation}}

\def\blu{\color{Blue}}
\def\mag{\color{Maroon}}
\def\red{\color{Red}}
\def\green{\color{ForestGreen}}
\def\prob{{\rm Prob}}

\def\reff#1{(\protect\ref{#1})}

\let\a=\alpha \let\b=\beta \let\ch=\chi \let\d=\delta \let\e=\varepsilon
\let\f=\varphi \let\g=\gamma \let\h=\eta    \let\k=\kappa \let\l=\lambda
\let\m=\mu \let\n=\nu \let\o=\omega    \let\p=\pi \let\ph=\varphi
\let\r=\rho \let\s=\sigma \let\t=\tau \let\th=\vartheta
\let\y=\upsilon \let\x=\xi \let\z=\zeta
\let\D=\Delta \let\F=\Phi \let\G=\Gamma \let\L=\Lambda \let\Th=\Theta
\let\O=\Omega \let\P=\Pi \let\Ps=\Psi \let\Si=\Sigma \let\X=\Xi
\let\Y=\Upsilon

\global\newcount\numsec\global\newcount\numfor
\gdef\profonditastruttura{\dp\strutbox}
\def\senondefinito#1{\expandafter\ifx\csname#1\endcsname\relax}
\def\SIA #1,#2,#3 {\senondefinito{#1#2}
\expandafter\xdef\csname #1#2\endcsname{#3} \else
\write16{???? il simbolo #2 e' gia' stato definito !!!!} \fi}
\def\etichetta(#1){(\veroparagrafo.\veraformula)
\SIA e,#1,(\veroparagrafo.\veraformula)
 \global\advance\numfor by 1
 \write16{ EQ \equ(#1) ha simbolo #1 }}
\def\etichettaa(#1){(A\veroparagrafo.\veraformula)
 \SIA e,#1,(A\veroparagrafo.\veraformula)
 \global\advance\numfor by 1\write16{ EQ \equ(#1) ha simbolo #1 }}
\def\BOZZA{\def\alato(##1){
 {\vtop to \profonditastruttura{\baselineskip
 \profonditastruttura\vss
 \rlap{\kern-\hsize\kern-1.2truecm{$\scriptstyle##1$}}}}}}
\def\alato(#1){}
\def\veroparagrafo{\number\numsec}\def\veraformula{\number\numfor}
\def\Eq(#1){\eqno{\etichetta(#1)\alato(#1)}}
\def\eq(#1){\etichetta(#1)\alato(#1)}
\def\Eqa(#1){\eqno{\etichettaa(#1)\alato(#1)}}
\def\eqa(#1){\etichettaa(#1)\alato(#1)}
\def\equ(#1){\senondefinito{e#1}$\clubsuit$#1\else\csname e#1\endcsname\fi}
\let\EQ=\Eq

\def\pp{{\bm p}}\def\pt{{\tilde{\bm p}}}


\def\\{\noindent}
\let\io=\infty
\def\ee{\end{equation}}
\def\be{\begin{equation}}

\def\VU{{\mathbb{V}}}
\def\EE{{\mathbb{E}}}
\def\UU{{\mathbb{U}}}
\def\N{\mathbb{N}}
\def\U{\mathbb{U}}
\def\GI{{\mathbb{G}}}
\def\TT{{\mathbb{T}}}
\def\C{\mathbb{C}}
\def\CC{{\mathcal C}}
\def\KK{{\mathcal K}}
\def\II{{\mathcal I}}
\def\LL{{\cal L}}
\def\RR{{\cal R}}
\def\SS{{\cal S}}
\def\NN{{\cal N}}
\def\HH{{\cal H}}
\def\GG{{\cal G}}
\def\PP{{\cal P}}
\def\AA{{\cal A}}
\def\BB{{\cal B}}
\def\FF{{\cal F}}
\def\v{\vskip.1cm}
\def\vv{\vskip.2cm}
\def\gt{{\tilde\g}}
\def\E{{\mathcal E} }
\def\EI{{\mathbb E} }
\def\I{{\rm I}}
\def\rfp{R^{*}}
\def\rd{R^{^{_{\rm D}}}}
\def\ffp{\varphi^{*}}
\def\ffpt{\widetilde\varphi^{*}}
\def\fd{\varphi^{^{_{\rm D}}}}
\def\fdt{\widetilde\varphi^{^{_{\rm D}}}}
\def\pfp{\Pi^{*}}
\def\pd{\Pi^{^{_{\rm D}}}}
\def\pbfp{\Pi^{*}}
\def\fbfp{{\bm\varphi}^{*}}
\def\fbd{{\bm\varphi}^{^{_{\rm D}}}}
\def\rfpt{{\widetilde R}^{*}}
\def\A{{{\mathcal O}}}
\def\ef{\mathfrak{f}}
\def\Ti{\mathfrak{T}}
\def\Mi{\mathfrak{M}}
\def\mm{\mathrm{ \mathbf{m}}}
\def\pp{\mathrm{ \mathbf{p}}}

\def\tende#1{\vtop{\ialign{##\crcr\rightarrowfill\crcr
              \noalign{\kern-1pt\nointerlineskip}
              \hskip3.pt${\scriptstyle #1}$\hskip3.pt\crcr}}}
\def\otto{{\kern-1.truept\leftarrow\kern-5.truept\to\kern-1.truept}}
\def\arm{{}}
\font\bigfnt=cmbx10 scaled\magstep1

\newcommand{\card}[1]{\left|#1\right|}
\newcommand{\und}[1]{\underline{#1}}
\def\1{\rlap{\mbox{\small\rm 1}}\kern.15em 1}
\def\ind#1{\1_{\{#1\}}}
\def\bydef{:=}
\def\defby{=:}
\def\buildd#1#2{\mathrel{\mathop{\kern 0pt#1}\limits_{#2}}}
\def\card#1{\left|#1\right|}
\def\proof{\noindent{\bf Proof. }}
\def\qed{ \square}
\def\trp{\mathbb{T}}
\def\trt{\mathcal{T}}
\def\Z{\mathbb{Z}}
\def\be{\begin{equation}}
\def\ee{\end{equation}}
\def\bea{\begin{eqnarray}}
\def\eea{\end{eqnarray}}
\def\kk{{\bf k}}
\def\Ti{\mathfrak{T}}
\def\Fi{\mathfrak{F}}
\def\Mi{\mathfrak{M}}
\def\begn{\begin{aligned}}
\def\egn{\end{aligned}}
\def\ti{{\rm\bf  t}}\def\mi{{\rm\bf m}}
\def\Va{{V^a_{\rm h.c.}}}
\def\Re{{\mathbb{R}}}
\def\T{{\mathcal{T}}}
\def\hL{{\L}}
\def\ev{\mathfrak{e}}
\def\obj{{\rm supp}}\def\fa{\FF}
\def\E{{\cal E}}
\def\EA{{E_\A}}
\def\0{\emptyset}
\def\Ni{\overline{\N}}

\title{On the zero-free region  for the chromatic polynomial of claw-free graphs with and without  induced square and induced diamond}
\author{
\\
Paula M. S. Fialho$^1$,  Aldo Procacci$^2$. \\
\\
\small{$^1$Departamento de Ci\^encia da Computa\c{c}\~ao UFMG, }
\small{30161-970 - Belo Horizonte - MG
Brazil}\\
\small{$^2$Departamento de Matem\'atica UFMG,}
\small{ 30161-970 - Belo Horizonte - MG
Brazil}\\
}

\maketitle

\begin{abstract}

Given a claw-free graph $\GI=(\VU,\EE)$ with maximum degree $\D$,
 we define the  parameter $\k\in [0,1]$ as $\k={\max_{v\in \VU}|I_v|\over \lfloor\D^2/4\rfloor}$ where $I_v$ is the set of all independent pairs in the neighborhood of $v$. We refer to $\k$  the pair independence ratio of $\GI$.
 We prove that for any claw-free graph $\GI$ with pair independence ratio  at most $\k$
the zeros of its chromatic polynomial $P_\GI(q)$  lie inside the  disk $D=\{q\in \C:~|q|< C_\k^0\Delta\}$,
where $C_\k^0$ is an increasing function of  $\k\in [0,1]$.
If $\GI$ is  also  square-free {and diamond free}, the function $C_\k^0$ can be replaced by a sharper function $C_\k^1$. These bounds constitute an improvement upon results
 recently given by Bencs and Regts \cite{BR}.

\end{abstract}

\section{Introduction and results}

\\Given $n\in \N$, we set $[n]=\{1,2,\dots,n\}$.
For any finite set $U$, we will
denote by $|U|$ its  cardinality and by $\mathcal{P}(U)$ the set of all subsets of $U$.  {Along this paper a  graph is an ordered  pair $\GI = (\VU,\EE)$ with vertex set $\VU$ finite, and edge set $\EE\subset \mathcal{P}_2(\VU)$ where $ \mathcal{P}_2(\VU)=\{e\subset \VU: |e|=2\}$
is the set of
all subsets of $\VU$ with cardinality 2.} Given a  graph $\GI=(\VU,\EE)$, a coloring  of $\VU$ with $q\in \N$ colors is a function $c: \VU\to [q]$. A coloring  of $\VU$ is called {\it proper} if for any edge $\{x,y\}\in \EE$ it holds that
$c(x)\neq c(y)$. Letting  $\KK^*_\VU(q)$ be the set of all proper colorings of $\VU$ with $q$ colors, the quantity $P_\mathbb{G}(q):= |\KK^*_\VU(q)|$, i.e. the number of proper colorings with $q$ colors of the graph $\mathbb{G}$, is, as a function of $q$, a polynomial called the \emph{chromatic polynomial} of $\GI$.

\\For a general finite graph $\GI$,
whose sole information is its maximum  degree $\D$,  Sokal \cite{sok01} proved that all the zeros of $P_\GI(q)$ lie
in the disk $|q| < C(\D)$, where $C(\D)$ is an explicitly computable function of $\D$
such that  $C(\D)/\D\le 7.963906$ for any $\D$.  Later,  Fern\'andez and Procacci \cite{FP2} improved the results  of Sokal by proving that
 $C(\D)/\D\le 6.9077$.

\\Recently  Jenssen, Patel and Regts \cite{JPR}, using a  representation of $P_\GI(q)$
known as  Whitney broken circuit Theorem \cite{W},
have improved the bound given in \cite{FP2} for the zero-free region of the chromatic polynomial. In \cite{JPR} the  authors show that for any
graph $\GI$ with maximum degree $\D$ and girth $g\ge 3$ the chromatic polynomial $P_\GI(q)$
is free of zero   outside the disk  $|q|<K_g \Delta$,
where $K_g$ is  a decreasing and explicitly computable function of $g$. In particular, the Jenssen-Patel-Regts result implies that
for any graph $\GI$ of maximum degree $\D$,
the chromatic polynomial is free of zeros if $|q|\ge  5.9315\D$. Then, Fialho, Juliano and Procacci \cite{FJP2} improved the results given in \cite{JPR}
by providing, for all pairs of integers $(\D,g)$  with  $\D\ge 3$ and $g\ge 3$, a positive number  $C(\D, g)$ such that
chromatic polynomial $P_\GI(q)$ of a graph $\GI$ with maximum degree $\D$ and finite girth $g$ is free of zero if $|q|\ge C(\D, g)$.
The bounds given in \cite{FJP2}
enlarge expressively the zero-free region in the complex plane of  $P_\GI(q)$  in comparison to  bounds in \cite{JPR},
for small values  of $\D$, while they coincide with those of \cite{JPR}  when $\D\to \infty$. Finally, very recently Bencs and Regts \cite{BR} have once again sensibly improved the bounds obtained  in \cite{JPR} and \cite{FJP2} lowering in  particular the zero-free region  of the chromatic polynomial of graph $\GI$  with maximum degree $\D$ by showing that all zeros of $P_\GI(q)$ lie  inside  the disk $|q|<4.25\D$.
Moreover in \cite{BR} the authors show that when $\GI$ is a claw-free graph then the zeros of $P_\GI(q)$ lie  in  the disk $|q|<3.81\D$.

\\In this note, we focus on the class of claw-free graphs and propose an improvement to the Bencs-Regts bound for the zero-free region of the chromatic polynomials. To state our main result we need to introduce some notations and definitions about claw-free graphs.

\\A graph $\GI=(\VU, \EE)$ is called {\it claw-free} if $K_{1,3}$ is not an induced subgraph
of $\GI$.

\\We recall that if $\GI=(\VU, \EE)$ is a claw free graph with maximum degree $\D$, then  for any $v \in \VU$, any independent subset of $\G_{\GI}(v)$ has size at most two. Therefore, the complement of $\G_{\GI}(v)$ must be triangle-free and, by Mantel's theorem \cite{Ma},  we get that { the number of independent pairs of vertices (a.k.a non-edges)  in $\G_{\GI}(v)$ satisfies}
\begin{equation}\label{eqMantel}
|\{S \subset \G_{\GI}(v);~~ |S|=2~ \mbox{and}~ S~ \mbox{independent in}~ \GI\}| \le~\left\lfloor\frac{\D^2}{4}\right\rfloor.
\end{equation}
We denote by $\GG_0$ the class of claw-free graphs. The class $\GG_0$ is {\it hereditary}, i.e every induced subgraph $H$ of a graph  $\GI \in \GG_0$ is also
claw-free ($H \in\GG_0$). We also
let $\GG_1 \subseteq \GG_0$ be the set of claw-free graphs that are also square-free and diamond-free graphs, namely  $\GI\in \GG_1$
if and only if $K_{1,3}$,  $C_4$ and $K_4\setminus e$ are not induced subgraphs of $\GI$. It is straightforward to verify
that $\GG_1$ is also a hereditary class, i.e if $\GI\in \GG_1$, then any induced subgraph H of $\GI$ is such that
$H\in \GG_1$.

\\Given a graph $\GI\in \GG_i$ with $i=0,1$, we introduce an other parameter of $\GI$ beyond the maximum degree $\D$, namely
the non negative number $\k$  defined as follows.

\begin{definition}\label{double}
Let $\GI=(\VU,\EE)$ be a graph and let $v\in \VU$. Let $\G_\GI(v)$ be the neighborhood of $v$
in $\GI$.
Let $I_v$  be the set of non-edges in  $\GI|_{\G_\GI(v)}$.
Then we define  the pair independence ratio of $\GI$ as
\be\label{kaa}
\k={\max_{v\in \VU}|I_v|\over \left\lfloor\frac{\D^2}{4}\right\rfloor}
\ee
\end{definition}
Observe that $\k$~is a number  taking values in the interval  $[0,1]$
and  the case $\k=0$ only happens if $\GI$ is the complete graph.

\begin{theorem}\label{teopri}
Let $\GI$ be a    graph with  maximum degree at most $\D\ge 3$ and pair independence parameter $\k$. If $\GI\in \GG_i$ with $i=0,1$, then  there exists a constant $C^i_{\k}$ such that $P_\GI(q) \neq 0$ for any $q \in \mathbb{C}$ satisfying
\[
|q|\ge C^i_{\k}\D
\]
where constant $C^i_{\k}$  increases with $\k$ and  ranges from $ C^1_{0}\le   3$ to  $C^0_{1}\le 3.81$.
\end{theorem}
\vskip.1cm


\begin{table}[t]
\centering
\caption{Constants $C^{i}_\k$ and corresponding minimizers $a_i^\ast$.}
\label{tab:deltafree_kappa_grid}
\begin{tabular}{|r|r|r|r|r|}
\hline
$\kappa$ & $C^{0}_\k$ & $C^{1}_\k$ & $a^{\ast}_{0}$ & $a^{\ast}_{1}$ \\
\hline
0.0 & 3.000000 & 3.000000 & 0.333333 & 0.333333 \\
\hline
0.1 & 3.169627 & 3.128158 & 0.355952 & 0.350761 \\
\hline
0.2 & 3.285039 & 3.214447 & 0.363300 & 0.356563 \\
\hline
0.3 & 3.377769 & 3.283304 & 0.367230 & 0.359765 \\
\hline
0.4 & 3.457121 & 3.341956 & 0.369749 & 0.361902 \\
\hline
0.5 & 3.527398 & 3.393730 & 0.371534 & 0.363490 \\
\hline
0.6 & 3.591011 & 3.440483 & 0.372885 & 0.364754 \\
\hline
0.7 & 3.649470 & 3.483371 & 0.373957 & 0.365812 \\
\hline
0.8 & 3.703793 & 3.523172 & 0.374839 & 0.366730 \\
\hline
0.9 & 3.754706 & 3.560437 & 0.375586 & 0.367548 \\
\hline
1.0 & 3.802747 & 3.595574 & 0.376232 & 0.368292 \\
\hline
\end{tabular}
\end{table}

\vskip.1cm
\\To prove Theorem \ref{teopri} we  use {an extension} of the Whitney broken circuit theorem proved in \cite{FJP}, namely a rewriting of the
 chromatic polynomial in term of forest  formed by  tree which are invariant under partition scheme and we use the Penrose partition scheme (see ahead for the definitions). We then  basically follow the scheme adopted in \cite{JPR,FJP2, BR}.

\\The class of (claw,diamond,$C_4$)-free graphs has been considered e.g. in \cite{LMZ} where it is proved that it has factorial speed. See also
\cite{dR} for more information on this class.
It is also worth mentioning that  (claw, $C_4$)-free graphs have been the subject of a recent study in \cite{CHKK} where the authors  prove a structure theorem for the class of (claw, $C_4$)-free graphs that are not quasi-line graphs.

\\In Section~\ref{secpre} we present definitions and tools that will be used throughout,
with emphasis on the Penrose partition scheme and its relation with respect to  the chromatic polynomial.
In Section \ref{sec3} we study the structure of Penrose trees in claw-free graphs and the corresponding generating-function bounds.
We finally prove Theorem~\ref{teopri}  in Section \ref{secproof}.

\numsec=2\numfor=1

\section{Preliminaries}\label{secpre}
\subsection{Trees, forests and partition schemes in a graph}

We start by establishing some definitions. A graph $\GI=(\VU, \EE)$
is {\it connected}
if for any pair $B, C$ of  subsets of $\VU$ such that
$B\cup C =\VU$ and $B\cap C =\emptyset$, there is at least an edge  $e\in \EE$ such
that $e\cap B\neq\emptyset$ and $e\cap C\neq\emptyset$. Given a vertex $u\in \VU$, we denote by $\G_{\GI}(u)=\{v\in \VU: \{u,v\}\in \EE\}$  the  neighborhood of $u$ in $\GI$. For any non-empty subset $E\subset \EE$, we define
$V_E=\{x\in \VU: x\in e ~{\rm for ~some~} e\in E\}$ and we say that such an $E \subset \EE$  is {\it connected }if the graph $\GI|_E=(V_E, E)$ is connected. Observe that any  non-empty $E\subset \EE$ can be written in a unique way as $E=\cup_{i=1}^kE_i$, for some integer $k\ge 1$, where $E_i$ is connected for all $i\in [k]$ and $V_{E_i}\cap V_{E_j}=\0$ for all $\{i,j\}\subset [k]$; the subsets $E_1, \dots, E_k$ are called the (maximal) {\it connected components} of $E$.

\\In this paper, a {\it tree} will be seen as a non empty connected subset of edges $\t\subset \EE$  such that $|V_\t|=|\t|+1$, and let $ \mathcal{T}_\GI$ be the set of all trees in $\GI$ plus  the empty set.
A tree $\t\in \mathcal{T}_\GI$ is said to be {\it rooted} if one of its vertices has been chosen as the  {\it root}.
If $x$ is a vertex of a rooted tree  $\t$, we let $d_\t(x)$ be the  {\it depth} of  $x$, i.e. the number of edges  in the path from the root to the vertex $x$, and we denote by $x^*$  the  {\it father} of $x$ in $\t$
 (i.e. $x^*$ is the vertex in the path from the root to  $x$ at depth $d_\t(x)-1$).

\\A {\it forest} in $\GI=(\VU,\EE)$ is a non-empty subset $F\subset \EE$ such that all its  connected components are trees and we   denote by $n_F\ge 1$ the number of trees in the forest $F$ (observe that if $n_F = 1$, then the forest $F$ is composed just by one tree $\t \in \mathcal{T}_\GI$). Let the symbol $\biguplus$ represents disjoint union. Then  each forest $F$ can be written in a unique way as a disjoint union of trees $F=\biguplus_{i=1}^{n_F}\t_i$ so that $V_F=\biguplus_{i=1}^{n_F}V_i$. Denote by $\mathcal{F}_\GI$ the set of all forests in $\GI$ plus the empty set, clearly we have that $\mathcal{T}_\GI\subset \mathcal{F}_\GI$.

\\A subset $E'\subset E$ is said a {\it spanning} subset of $E \subset \EE$ if $V_{E'}=V_{E}$.
Given a subset $R\subset \VU$ such that $|R| \ge 2$, we define $\EE|_R=\{e\in \EE: |e\cap R|=2\}$ and represent by $\mathcal{C}^{\rm sp}_{\EE|_R}$ the set of all connected spanning subsets of ${\EE|_R}$. Analogously, we define
$\mathcal{T}^{\rm sp}_{\EE|_R}\subset {\mathcal{C}}^{\rm sp}_{\EE|_R}$ as the set of  all  spanning trees in $\EE|_R$. Finally, if $E\subset E'\subset \EE$, we set   $$[E,E']=\{E''\subset \EE: E\subseteq E''\subseteq E'].$$
It is  convenient  to organize the collection of connected spanning subsets of a graph in a systematic way by giving  a so-called {\it partition scheme}, whose definition we give below.

{\begin{definition}\label{defia}
A partition scheme in $\GI$ is a map  $\mm: \mathcal{F}_\GI\to \mathcal{P}(\EE)$
 such that
\vv
(1) $\mm(\0)=\0$.
\vv
(2) For each $\t\in {\mathcal{T}}_\GI$,  $\t\subset \mm(\t)$ and $V_\t=V_{\mm(\t)}$.
\vv
(3) For each $R\subset \VU$ with $|R| \ge 2$,  $\mathcal{C}^{\rm sp}_{\EE|_R}=\biguplus_{\tau\in \mathcal{T}^{\rm sp}_{\EE|_R}}[\tau, \mm(\tau)]$.
\vv
(4) For any  forest  $F\in \mathcal{F}_\GI$ such that $F={\biguplus_{i=1}^{ n_F}\t_i}$, $\mm(F)= \bigcup_{i=1}^{n_F}\mm(\t_i)$.
\end{definition}}

\\Observe that by item $\it (3)$ a partition scheme decomposes the family of connected spanning edge sets into disjoint intervals indexed by spanning trees and  item $\it (4)$ extends this notion to forests.

\\Recently in \cite{FJP} the authors, inspired by Whitney's broken circuit free Theorem (see \cite{W}), showed that any partitions scheme can be used to define the chromatic polynomial of a graph, by proving the following theorem.

\begin{theorem}\label{partschem} Let $\GI=(\VU,\EE)$ be a graph and let  $P_\GI(q)$ be the chromatic polynomial of $\GI$. Let $\mathrm{\mathbf{m}}$ be a partition scheme in $\GI$  according to Definition \ref{defia}.
 Then
\be\label{pgen}
P_\GI(q)=q^{{{|\VU|}}} \sum_{F\in \mathcal{F}^\mm_\GI} \left(-{1\over q}\right)^{|F|}
\ee
where $\mathcal{F}^{\mm}_\GI$ is the subset of $\mathcal{F}_\GI$ formed by  all  forests $F$ in $\GI$ such that $\mm(F)=F$.
\end{theorem}

\\The representation \reff{pgen} for $P_\GI(q)$ given in Theorem \ref{partschem} will be used in the proof of Theorem \ref{teopri}, together with  the so-called {\it Penrose
partition scheme}, {originally proposed by Oliver Penrose in \cite{pen67}}.

{\begin{definition}[Penrose  partition scheme]\label{pens}
Given a graph $\GI=(\VU,\EE)$, fix a total ordering $\succ$ in $\VU$ in such a way that any tree $\t\in \mathcal{T}_\GI$ is considered to be rooted at its least vertex. Let $\bm\pp: {\mathcal{T}}_\GI\to \mathcal{P}(\EE): \t\mapsto \pp(\t)$ be the map such that
${\pp}( \tau )=\t\cup E_\tau^{\bm p}$, where
$$E_\tau^{\pp}=\{\{x,y\}\subset \EE|_{V_\t}\setminus \t; ~~\mbox{either}~~ d_\t(x)=d_\t(y), ~~\mbox{or}~~ d_\t(y)=d_\t(x)-1~~\mbox{and}~~ y\succ x^*\}.$$
Then, the Penrose  partition scheme is the map
$\pp:\mathcal{F}_\GI\to\mathcal{P}(\EE)$ such that $\pp(F)= \biguplus_{i=1}^{n_F}\pp(\t_i)$ for any non-empty
$F\in \mathcal{F}_\GI$ such that  $F=\biguplus_{i=1}^{n_F}\t_i$.
\end{definition}}

\\{A detailed  proof that the map $\bm p$ is indeed a partition scheme according to Definition \ref{pens} can be found in \cite{FP}.}

\\A tree $\t\in \mathcal{T}_\GI$ such that
 $\pp(\t)=\t$ will be called a {\it Penrose tree} and a forest $F\in \mathcal{F}_\GI$ such that $\pp(F)=F$ is called a Penrose forest. We shortly denote by
 $ \mathfrak{T}_{\GI}=\{\t\in \mathcal{T}_\GI: \pp(\t)=\t\}\cup\0$ and
 $ \mathfrak{F}_{\GI}=\{F\in \mathcal{F}_\GI: \pp(F)=F\}\cup\0$,   the set of all Penrose trees (empty set included) and Penrose forests (empty forest included) in $\GI$ respectively.
Therefore, for any $z\in \C$  and any graph $\GI$, according to Theorem \ref{partschem}, we have that
\[
P_\GI(q)=q^{|\VU|}{F}_\GI(-1/q)
\]
where
\be\label{FGz}
{F}_\GI(z) =\sum_{F \in \mathfrak{F}_{\GI}}z^{|F|}.
\ee
Observe that $P_\GI(q) \neq0$ as soon as the function defined in \reff{FGz} is not zero. {Therefore}, instead of studying the chromatic polynomial directly, { we can investigate the zero-free region}  of
$F_\GI(z)$ in order to prove Theorem \ref{teopri}.

\\Let us now give some   notations and properties regarding $F_\GI(z)$.

\\Given the graph $\GI=(\VU,\EE)$, we abuse notation and also write ${F}_\VU(z)$ in place of $F_\GI(z)$, when we need to emphasize the dependency on the vertex set, as well as ${F}_U(z)$ in place of ${F}_{\GI|_U}(z)$, for any $U\subset \VU$.
Given  $S \subseteq \VU$, let $\mathfrak{F}_{\GI,S} \subset \mathfrak{F}_{\GI}$ be the set of all Penrose forests $F$ of $\GI$  such that every tree  of $F$ contains at least one vertex of $S$ (the empty set is included in  $\mathfrak{F}_{\GI,S}$). Moreover, for any $U \subseteq \VU$ and $S\subset U$ we will write shortly $\mathfrak{F}_{U}$,
 $\mathfrak{F}_{U,S}$   instead of $\mathfrak{F}_{\GI|_U}$,  $\mathfrak{F}_{\GI|_U,S}$  respectively. Note that if $S=\{v\} \subset U$ then $\mathfrak{F}_{U,u}$ is the set of all Penrose trees in $\GI|_U$ containing $v$, which we denote simply by $\mathfrak{T}_{U,v}$ (observe that the empty set is contained in $\mathfrak{T}_{U,v}$).

 \\Given $U\subseteq \VU$,  observe that, for any non empty set $S\subset U$, we have
\be\label{FU}
{F}_U(z)= \sum_{F\in \mathfrak{F}_{U, S}}z^{|F|}{{F}_{U\setminus (S\cup V_F)}}(z),
\ee
hence for all $z\in \C$ such that ${F}_U(z)\neq 0$, we can write
\be\label{idun}
\sum_{F\in \mathfrak{F}_{U,S}}z^{|F|}{{F}_{U\setminus (S\cup V_F)}(z)\over {F}_U(z)}=1, ~~~~\mbox{for all $U\subseteq \VU$ and $S\subset U$}.
\ee
In particular, given any  vertex $u\in \VU$,  identity  \reff{FU} specialized to $U=\VU$ and $S=\{u\}$ can be written as
\begin{equation}
\label{eq:recur}
{F}_\VU(z) = {F}_{\VU-u}(z) + \sum_{\t \in \mathfrak{T}_{\VU,u}\atop |\t|\ge 1}z^{|\t|}{F}_{\VU-u -V_\t}(z)\, .
\end{equation}
with $\mathfrak{T}_{\VU,u}$ the set of all Penrose trees of $\GI$ containing $u$, as defined above. If $u\in \VU$ and $z\in \C$ are such that  ${F}_{\VU-u}(z)\neq 0$, we define
\be\label{rug}
R^u_\GI(z)= \frac{{F}_{\VU}(z)}{{F}_{\VU-u}(z)} - 1
\ee
and \reff{eq:recur} implies that
\begin{equation}
\label{RuG}
R^u_\GI(z)= \sum_{\substack{\t \in \mathfrak{T}_{\VU,u} \\ |\t| \geq 1 }}z^{|\t|}
\frac{{F}_{\VU-u-V_\t}(z)}{{F}_{\VU-u}(z)}.
\ee

\\The relations \reff{FU}-\reff{RuG} will be used to find a disk where the  $F_\GI(z)$ is different from zero. In next  section we will explore properties of Penrose trees.

\\We  conclude this section by showing  a property regarding Penrose trees  which will be used ahead toward the proof of Theorem \ref{teopri}.

\begin{lemma}\label{lem:case2-obstruction}
Consider a graph $\GI=(\VU,\EE)$ and an ordering $\prec$ in $\VU$ such that $u$ is the least vertex. Let
$S=\{v_1,v_2\}\subseteq \Gamma_\GI(u)$ be an independent pair with {$v_1 \prec v_2\prec w$ for all $w\in \VU-u$}.  Consider a forest { $F \in \mathfrak{F}_{\VU-u}$ formed  by two trees $\t_1$ and $\t_2$ with
$\t_i\in  \mathfrak{T}_{\VU-u,v_i}$ for $i=1,2$}.
Then $F\cup uS$ is \emph{not} a Penrose tree in $\GI$ if and only if at least one of the following holds:
\vv
\\a) There is an edge $\{x,y\}\in \EE$ with $x\in V_{\tau_1}$ and $y\in V_{\tau_2}$ such that $d_{F\cup uS}(x)=d_{F\cup uS}(y)$.
\vv
\\b) There is an edge $\{x,y\}\in \EE$ with $x\in V_{\tau_1}$ and $y\in V_{\tau_2}$ such that either $d_{F\cup uS}(x)-1=d_{F\cup uS}(y)$ and $y\prec x^*$ or $d_{F\cup uS}(y)-1=d_{F\cup uS}(x)$ and $x\prec y^*$.
\vv
\\c) There is an edge $\{v_2, x\}\in \EE$ where $x$ is a \emph{child} of $v_1$ in $\tau_1$.
\vv
\\In particular, if neither a) nor b) nor c) occur, then $F\cup uS$ is a Penrose tree in $\GI$.
\end{lemma}

\\{\bf Proof.} By Definition~\ref{pens}, the Penrose map for a rooted tree $\t$ adds exactly two types of edges in $\EE \setminus \t$:

\\(i) edge $\{v,w\}$, if $v$ and $w$ are vertices in $\t$ such that $d_\t(v)=d\t(w)$,

\\(ii) edge $\{v',w\} \in \EE$, if $\{v, v'\}\in V_\t$ and $v$ and $w$ are vertices in $\t$ such that $d_\t(v)=d\t(w)=d_\t(v')-1$.

\\Since each $\tau_i$ { (with root $v_i$ by the ordering chosen)} is already Penrose, any violating edge in $\pp(F\cup uS)$ must connect the two branches below $u$,
and these are precisely \textnormal{(a)},   \textnormal{(b)} and \textnormal{(c)}. If neither occurs,  hence $\pp(F\cup uS)=F\cup uS$ and the union is Penrose.

\hfill$\Box$

\begin{corollary}\label{cor:one-trivial}
Let $\GI\in\mathcal G_1$ and fix a total an ordering $\prec$ in $\VU$ such that $u$ is the least vertex. Let
$S=\{v_1,v_2\}\subseteq \Gamma_G(u)$ be an independent pair with { $v_1 \prec v_2\prec w$ for all $w\in \VU-u$}. If $\tau_1$ is any Penrose tree rooted at $v_1$ and $\tau_2=\emptyset$ (rooted at $v_2$), then
$F\cup uS$ is a Penrose tree in $\GI$.
\end{corollary}

\\{\bf Proof.}
By Lemma~\ref{lem:case2-obstruction}, with one branch trivial the only possible violation is an edge
$\{v_2,w\}$ where $w$ is a child of $v_1$ in $\tau_1$. However  such an edge cannot exist in a  ${(C_4,K_4-e)}$-free graph. Indeed,  if $w$ is not a neighbor of $u$ the existence of the edge $\{w,v_2\}$ would imply that $\GI$ would contain the square with edges $\{u,v_1\},\{v_1,w\},\{w,v_2\},\{v_2,u\}$. On the other hand if $w$ is a neighbor of $u$ then
existence of the edge $\{w,v_2\}$ would imply that $\GI$ would contain the diamond  with edges $\{u,v_1\},\{v_1,w\},\{w,v_2\},\{v_2,u\}, \{w,u\}$.
Hence $\pp(F\cup uS)=F\cup uS$.

\hfill$\Box$

\\Let $\GI=(\VU,\EE)\in \GG_i$ ($i=0,1$) with maximum degree $\D\ge 2$. For $y\ge 0$ we define the {\it generating function of Penrose trees} in $\GI$ rooted at $v$, as follows:
\be\label{tgv}
T_{\GI,v}(y):=\sum_{\substack{\t\in \mathfrak{T}_{\GI,v}}} y^{|\t|},
\ee
where recall that $\mathfrak{T}_{\GI,v}$ is the set of all Penrose trees in
 $\mathcal{T}_\GI$  that contains $v$ and $|\t|$ is the number of edges in $\t$. Note that the sum in the r.h.s. of \reff{tgv} is a polynomial in $y$ starting with    1, which occurs when $\t=\0$.

\section{A bound on $T_{\GI,v}(y)$ when $\GI$ is a claw-free graph}\label{sec3}

\\Fix an integer $d \geq 2$, and let $U(d)$ denote the infinite rooted tree in which each vertex has exactly $d$ children and let $\UU$ denote the set of vertices of $U(d)$. Now, for each vertex $v\in U(d)$, let $C_v\subset \UU$ be the set of vertices which are children of $v$ in $U(d)$  and let
$$
P_v=\{\{w,w'\}\subset \UU:~ w\in C_v ~{\rm and}~ w'\in C_v\}
$$
Note that $|P_v|= {d\choose 2}$.

\\Suppose  that for each $v\in U(d)$ it is fixed a set
$I_v\subset P_v$ with cardinality $|I_v|= m$ where $m=\left\lfloor{d^2\over 4}\right\rfloor$ and let $u_n(d,m)$ be the number of finite  subtrees $\tau \subseteq U(d)$ such that:

\\-  {the root of $U(d)$ belongs to $V_\t$}

\\- $|V_\t|=n$;

\\- {
any vertex $v\in V_\t$ has at most two children};

\\- {if $v\in V_\t$ has two children, say  $v_1$ and $v_2$, then the pair $\{v_1,v_2\}$ belongs to $I_v$.}

\\For $y\ge 0$, we define
\be\label{ux}
u^m_{d}(y):=\sum_{n\ge 1} u_{n}(d,m) y^{n}
\ee
and observe that\be\label{Us}
u^m_{d}(y)=yZ^m_\D(u^m_{d}(y)),
\ee
where
\begin{equation}\label{ZD}
Z^m_d(u)=
1+du+m u^2
\end{equation}
\vv

\\we also define
\be\label{W(x)}
w_d^m(y):= \sum_{n\ge 1}  u_{n}(d,m) y^{n-1}= {u_d^m(y)\over y}= Z^m_{d}(u_d^m(y)).
\ee

\begin{proposition}\label{A1}
For any $d\geq 2$ and any $m\ge 0$, the series \eqref{ux} converges for all $y\in [0,{1\over 2\sqrt{m}+d})$ and
\be\label{wws}
w_d^m(y)=\displaystyle{2\over  1-dy+\sqrt{(1-dy)^2-4m y^2}}.
\ee
\end{proposition}
\\{\bf Proof}.
If $m=0$, then no vertex of $\tau$ can have two children (since that would require the existence of an
independent pair), hence every vertex has at most one child and any
$\tau\subseteq U(d)$ counted by $u_n(d,0)$ is necessarily a rooted path with length $n-1$. Therefore,
\[
u_d^0(y)
=\sum_{n\ge 1} u_n(d,0)\,y^n
=\sum_{n\ge 1} d^{\,n-1}y^n
=\frac{y}{1-dy},
\qquad \text{for } 0\le y < \frac{1}{d}.
\]
Moreover, \(
w_d^0(y)
=\sum_{n\ge 1} d^{\,n-1}y^{n-1}
=\frac{1}{1-dy}\), for \( 0\le y < \frac{1}{d}.
\)


\\We are left with  the case $1\le m\le \lfloor{(\D-1)^2\over 4}\rfloor$.
The relation \eqref{Us} implies that
$$
u_d^m(y)=q^{-1}(y)
$$
where
\begin{equation*}
q(u)\;=\;{u \over Z_d^m(u)}.
\end{equation*}
For $u\ge 0$,  the function $q$ starts at $0$ and it is  strictly increasing in the interval
$[0,{1\over \sqrt{m}}]$, at the end of it reaches its maximum value $q({1\over \sqrt{m}})={1\over 2\sqrt{m}+d}$.
Therefore, $q$ is a bijection from $[0,{1\over \sqrt{m}})$ onto $[0,{1\over 2\sqrt{m}+d})$ and so $u_d^m(y)$ satisfying (\ref{Us}) is defined
in the whole interval $y\in [0,{1\over 2\sqrt{m}+d})$.
That being so,   the series (\ref{ux}) converges in the interval  $[0,{1\over 2\sqrt{m}+d})$.

\\To show \reff{wws}, from $u=yw$ and $w=1+du+m u^2$ we get that $w=1+dyw+m(yw)^2$. Whence, solving the quadric  choosing the branch such that $w(0)=1$ we get
$$
w_d^m(y)={1-dy- \sqrt{(d^2-4m)y^2-2dy+1}\over 2m y^2}
={2\over  1-dy+\sqrt{(1-dy)^2-4m y^2}}
$$
ending the proof
of  the proposition.

~~~~~~~~~~~~~~~~~~~~~~~~~~~~~~~~~~~~~~~~~~~~~~~~~~~~~~~~~~~~~~~~~~~~~~~~~~~~~~~~~~~~~~~~~~~~~~~~~~~~~~~~~~~~~~~~~~~$\Box$

\\Let now $\GI=(\VU,\EE)$ be a claw-free graph with maximum degree $\D$.
Then, for any $v\in \VU$ and any $w\in \G_\GI(v)$, by Mantel's Theorem the number of
independent pairs in $\G_\GI(v)\setminus\{w\}$ is  less than or equal to $\lfloor(\D-1)^2/4\rfloor$.  Let now $u$ be a vertex of $\GI$
{ with degree at most $\D-1$} and consider a tree
$\t\in \mathfrak{T}_{\GI,u}$.
According to definitions \reff{ux}-\reff{W(x)}, we can bound $T_{\GI,v}(x)$ as
\be\label{tvt}
T_{\GI,v}(y)\le  w_{\D-1}^{\lfloor(\D-1)^2/4\rfloor}(y).
\ee
when $y$ varies in the interval $[0, {1\over 2\sqrt{\lfloor(\D-1)^2/4}\rfloor+\D-1})$.
 \begin{lemma}\label{le2} Given a  claw-free graph  with maximum  degree at most $\D\ge 3$, let $v\in \VU$ be a vertex of $\GI$
with degree at most $\D-1$.
Let, for all $y\in [0,{1\over 2(\D-1)})$,
\be\label{gide}
g_{\D}(y)={4\over  (1+\sqrt{1-2(\D-1)y})^2}
\ee
Then,  for all $y\in [0, {1\over 2(\D-1)}]$ and all $\D\ge 3$,
\be\label{ineqv}
T_{\GI,v}(y)\le g_{\D}(y)
\ee
\end{lemma}

\\{\bf Proof}. Proposition \ref{A1} and   inequality \reff{tvt} immediately imply that,
for all $y\in [0, {1\over 2(\D-1)}]$,
\be\label{ineqv}
\begin{aligned}
T_{\GI,v}(y)& \le \displaystyle{2\over  1-(\D-1)y+\sqrt{[1-(\D-1)y]^2-4\lfloor(\D-1)^2/4\rfloor y^2}}\\
&\le \displaystyle{2\over  1-(\D-1)y+\sqrt{[1-(\D-1)y]^2-(\D-1)^2 y^2}}\\
&=\displaystyle{2\over  1-(\D-1)y+\sqrt{1-2(\D-1)y}}\\
&={4\over  (1+\sqrt{1-2(\D-1)y})^2}
\end{aligned}
\ee

\vv

\section{Proof of Theorem \ref{teopri}}\label{secproof}

We recall that in Section \ref{secpre} we showed that
\(
P_\GI(q)=q^{|\VU|}{F}_\GI(-1/q)
\)
where
\be
{F}_\GI(z) =\sum_{F \in \mathfrak{F}_{\GI}}z^{|F|}.
\ee
Therefore $P_\GI(q) \neq0$ as soon as ${F}_\GI(-1/q) \neq 0$. We will use relations \reff{FU}-\reff{RuG} to find a disk where the Penrose expansion $F_\GI$ is different from zero.

\subsection{A key proposition and an auxiliary lemma}
\\Let us start by giving some definitions.

\\We set
\be\label{ginf}
h(x)= {4\over  (1+\sqrt{1-2x})^2}
\ee
and observe that for any $x\in [0,{1\over 2}]$ and for any $\D\ge 2$ we have that
\be\label{ginf2}
g_\D({x\over \D})\le h(x)
\ee

\\Given   $a \in (0,1)$,  $\D\ge 3$,  $x\in [0, {1\over 2}]$,  $i\in\{0,1\}$ and  $\k\in [0,1]$, we define

\be\label{kdax}
K^i_{\k}(a,x)=
(1-a) x+\k {x^2\over 4}\left[(1-a)^2+\Big(h\left({x}\right)-1\Big)\Big(h\left({x}\right)-i\Big)\right]
\ee
\be\label{xda}
x^i_{\k}(a)=\sup\{x\in [0,{1/2}]: ~K^i_{\k}(a,x)\le a\}
\ee
\be
C_{\k}^i(a)= {1\over (1-a)x^i_{\k}(a)}
\ee
and
\be\label{zda}
z^i_{\k}(a)=
{1\over C_{\k}^i(a)\D}
\ee

\vv

\noindent\textbf{Remark.}
%
%
Note that the function $K^i_{\k}(a,x)$ is continuous and strictly increasing in $x$. Moreover
for any fixed $a\in (0,1)$ we have,
for  $\k<\k'$
$$
K^i_{\k'}(a,x)\ge K^i_{\k}(a,x), ~~~\forall x\in \left[0,{1\over 2}\right].
$$
Therefore, according to \reff{xda}, for each fixed $a\in(0,1)$ and for all $\D\geq 3$, if $0\le\k<\k'\le 1$, we have that $x^i_{\k'}(a)\;\leq\; x^i_{\k}(a)$
which in turn implies  that $C^i_{\k}(a)\leq C^i_{\k'}(a)$ for all $a\in (0,1)$.

 \begin{proposition}\label{pro1}
Let $\GI\in \GG_i$ for $i=0,1$, with maximum degree at most $\D\ge 3$ and pair independence ratio at most $\k$.
For every $a\in(0,1)$ and any $u\in \VU$, if $|z| \leq  z^i_{\k}(a)$, with $ z^i_{\k}(a)$ defined in \reff{zda},
then
\be\label{claim1}
{F}_\GI(z)\neq 0
\ee
and
\be\label{rlea}
|R^u_\GI(z)|\le a.
\ee
\end{proposition}

\\Theorem \ref{teopri} follows straightforwardly from Proposition \ref{pro1}.
Indeed, Proposition \ref{pro1} implies that for any graph $\GI\in\GG_i$ (for $i=0,1$)  with maximum degree $\D\ge 3$
 we have that ${F}_\GI(z)\neq 0$  whenever  $z\in \C$  is such that $|z| \leq z^i_{\k}(a)$
 for  any   $a\in (0,1)$. Then,
 since $P_\GI(q)= q^{|\VU|}{F}_\GI(-1/q)$, this implies  that,  for  any   $a\in (0,1)$,
the chromatic polynomial $P_\GI(q)$ is
free of zeros as soon as
 $q\in \C$ satisfying
 $$|q| \ge C^i_{\k}(a)\D.$$
Setting
$$
C^i_{\k}=\inf_{a\in (0,1)}C^i_{\k}(a)
$$
we conclude that
$P_\GI(q)\neq 0$ as soon as $|q| \geq C^i_{\k}\D$.

\\Values of $C^i_{\k}$ can be easily evaluated via computer-assisted computation; a selection of these values is listed in Table 1."

\\It remains to prove Proposition \ref{pro1}. Let us first state and prove
the following lemma.
\begin{lemma}\label{lem37}
Let $\GI=(\VU,\EE)$ be a  graph   with maximum degree at most $\D\ge 3$ and pair independence ratio at most $\k$ such that $\GI\in \GG_i$ for $i=0,1$.
Given $u\in \VU$,   let  $\GI'=(\VU-u, \EE|_{\VU-u})$ and  $\VU'=\VU-u$.
Assume that Proposition \ref{pro1} is true for the  graph $\GI'$ and all its subgraphs, then, given any $a\in (0,1)$ and
 any  $A\subset \VU'$,
\[
\left| \frac{{F}_{\VU'-A}(z)}{{F}_{\VU'}(z)} \right|\leq (1-a)^{-|A|}~~~~~\mbox{ for any $|z|\le z^i_{\k}(a)$}.
\]
\end{lemma}
\\{\bf Proof}.
Let  $A = \{v_1, \ldots, v_k\}  \subseteq \VU'$ and for $j\in [k]$, let $A_j = \{v_1, \ldots, v_j \}$ and $A_0 = \emptyset$.
Then by  hypothesis, for any $j\in [k]$ and for any  $|z|\le  z^i_{\k}(a)$, we have that ${F}_{\VU'-A_{j}}(z)\neq 0$ and
\be\label{indu2}
\left| \frac{{F}_{\VU'-A_{j}}(z)}{{F}_{\VU'-A_{j+1}}(z)} - 1 \right|\le a~~\Longrightarrow~~
\left|\frac{{F}_{\VU'-A_{j+1}}(z)}{{F}_{\VU'-A_{j}}(z)} \right|\le {1\over 1-a}~.
\ee
Hence
\[
\left| \frac{{F}_{\VU'-A}(z)}{{F}_{\VU'}(z)} \right|
= \prod_{i=1}^{k} \left| \frac{{F}_{\VU'-A_{i}}(z)}{{F}_{\VU'-A_{i-1}}(z)} \right|
\leq (1-a)^{-|A|},~~ ~~~~ \forall z\in  \C: ~ |z|\le z^i_{\k}(a).
\]
 $~~~~~~~~~~~~~~~~~~~~~~~~~~~~~~~~~~~~~~~~~~~~~~~~~~~~~~~~~~~~~~~~~~~~~~~~~~~~~~~~~~~~~~~~~~~~~~~~~~~~~~~~~~~~~~~~~~~~~~~~~~~\Box$

\subsection{Proof of Proposition \ref{pro1}}

\\Let $\GI=(\VU,\EE)$ be  a graph belonging to the class $\GG_i$ (with $i=0,1$) with maximum degree at most $\D$ and pair independence ratio at most $\k$. The proof will be done by induction on the number of vertices.
The induction is justified by the hereditariness of $\GG_0$ and $\GG_1$, together with the fact that the properties of maximum degree at most $\D$ and pair independence ratio at most $\k$ are inherited by subgraphs.
\vv
\\{\it If $|\VU|=2$:} then  $F_\VU(z)$ is either $1$ (if $\GI$ has no edge) or $(1+z)$ (if $\GI$ is a single edge) and thus $R_\GI^u(z)=0$ or  $R_\GI^u(z) = z$.
If $F_\VU(z)=1$   and
$R_\GI^u(z)=0$ then $F_\VU(z)\neq 0$ and $|R_\GI^u(z)|=0$ for  all $z\in \C$. If  $F_\VU(z)=1+z$   and $R_\GI^u(z)=z$ then  $F_\VU(z)\neq 0$ for  all $|z|< 1$. In both case Proposition \ref{pro1} is true since, by construction,  $z_{\k}^i(a)<1$ for any $a\in (0,1)$, $\D\ge 3$, any
$\k\in[0,1]$ and any $i\in \{0,1\}$.

\vv
\\{\it If $|\VU|\ge 3$:} given $u\in \VU$,   let us set  $\GI'=(\VU-u, \EE|_{\VU-u})$ and  $\VU'=\VU-u$. Then by the induction hypothesis, given $a\in (0,1)$, we have that $F_{\VU'}(z)\neq 0$ for $|z|\le z^i_{\k}(a)$.
In light of this, to conclude the proof  of Proposition \ref{pro1},  it is sufficient to prove    that for  any $a\in (0,1)$
\be\label{aux}
|R^u_\GI(z)|\le a~~~~~~~\mbox{for all $z\in \C$ such that  $|z|\leq z^i_{\k}(a)$}.
\ee
Indeed, if \reff{aux} holds, by the induction hypothesis  we have that ${F}_{\VU-u}(z)\neq 0$ for $|z|\leq z^i_{\k}(a)$, so that
$$
|R^u_\GI(z)|\le a~~\Longrightarrow~~\Big|\frac{{F}_\VU(z)}{{F}_{\VU-u}(z)} - 1 \Big|\le
a~~\Longrightarrow~~\Big|\frac{{F}_\VU(z)}{{F}_{\VU-u}(z)}\Big|\ge 1-a > 0
~~\Longrightarrow~~{F}_\VU(z)\neq 0.
$$
\\We now prove inequality  \reff{aux}. Let us fix an ordering of the set $\VU$ in such a way  that $u$ is the least vertex and the vertices  in
$\G_{\GI}(u)$
are  the least in the set $\VU-u$.
We stress that
such choice of the ordering  in $\VU$ is without loss of generality  because for any $U\subset \VU$ the function $F_U(z)$ does not depend on the order established in $\VU$.

\\For any $S\subset \G_\GI(u)$,
let $uS= \{e=\{u,s\}\in \EE: s \in S\}$. Moreover,
let $\mathcal{F}_{\VU',S}$ be the set of all forests of $\GI-u$ (either Penrose or not, empty set is included)  such that each tree  of $F$ contains  at least a vertex of $S$. Then
 the expression of  $R^u_\GI(z)$ given in \reff{RuG} can be written as follows.
\be
\label{Rugiz}
\begin{aligned}
R^u_\GI(z)
& = \sum_{\substack{\t \in \mathfrak{T}_{\GI,u} \\ |\t| \geq 1 }}z^{|\t|}\frac{{F}_{\VU'-V_\t}(z)}{{F}_{\VU'}(z)}\\
& =\sum_{S \subseteq \G_\GI(u) \atop S \not= \emptyset} z^{|S|}
\sum_{F \in \mathcal{F}_{\VU',S} \atop  F \cup uS\in\Ti_{\GI,u}}
z^{|F|}\frac{F_{\VU'\setminus(S\cup V_F)}(z)}{{F}_{\VU'}(z)}\\
\end{aligned}
\ee
where  we recall that $\mathfrak{T}_{\GI,u}$ is the set of all Penrose trees contained  in $\GI$ and containing $u$ (empty set included).

\\Now observe that,  due to our choice of the order in $ \VU$ we have
$$
\sum_{F \in \mathcal{F}_{\VU',S} \atop  F \cup uS\in\Ti_{\GI,u}}
z^{|F|}\frac{F_{\VU'\setminus(S\cup V_F)}(z)}{{F}_{\VU'}(z)}=
\begin{cases}
\sum_{F \in \mathfrak{F}_{\VU',S} \atop  F \cup uS\in\Ti_{\GI,u}}
z^{|F|}\frac{F_{\VU'\setminus(S\cup V_F)}(z)}{{F}_{\VU'}(z)} & {\rm if~} S ~{\rm is ~independent}\\
0 & {\rm otherwise}
\end{cases}
$$
{Indeed if $S$ is not independent, any tree $\t$ which is the union of a forest $F\in \mathcal{F}_{\VU',S}$ and $uS$ is not Penrose.
On the other hand, if $S$ is independent  but  $F\in  \mathcal{F}_{\VU',S}$ is not Penrose then  also
$F\cup\{uS\}$  is not Penrose.} So we can write
$$
R^u_\GI(z)
=\sum_{S \subseteq \G_\GI(u) \atop S \not= \emptyset, \text{ind}} z^{|S|}
\sum_{F \in \Fi_{\VU',S} \atop  F \cup uS\in\Ti_{\GI,u}}
z^{|F|}\frac{F_{\VU'\setminus(S\cup V_F)}(z)}{{F}_{\VU'}(z)}\,.
$$
\\Recalling that in a claw-free graph $\GI=(\VU,\EE)$, for any $u\in \VU$, any subset $S\subseteq \G_\GI(u)$ with $|S|\ge 3$ cannot be independent,
the above sum counts only when $|S|\in \{1,2\}$. If $S=\{v\}$ then, by \reff{idun},
$$
\sum_{F \in \Fi_{\VU',\{v\}} \atop  F \cup \{u,v\}\in\Ti_{\GI,u}}
z^{|F|}\frac{F_{\VU'\setminus(\{v\}\cup V_F)}(z)}{{F}_{\VU'}(z)} =
\sum_{\t \in \Ti_{\VU',\{v\}} }
z^{|\t|}\frac{F_{\VU'\setminus(\{v\}\cup V_\t)}(z)}{{F}_{\VU'}(z)} =1
$$
and, if $|S|= 2$, then
$$
\begin{aligned}
\sum_{F \in \Fi_{\VU',S} \atop  F \cup uS\in\Ti_{\GI,u}}
z^{|F|}\frac{F_{\VU'\setminus(S\cup V_F)}(z)}{{F}_{\VU'}(z)}
&=\sum_{F \in \Fi_{\VU',S} }
z^{|F|}\frac{F_{\VU'\setminus(S\cup V_F)}(z)}{{F}_{\VU'}(z)} -
\sum_{F \in \Fi_{\VU',S} \atop  F \cup uS \not\in\Ti_{\GI,u}}
z^{|F|}\frac{F_{\VU'\setminus(S\cup V_F)}(z)}{{F}_{\VU'}(z)}\\
& = 1-\sum_{F \in \Fi_{\VU',S} \atop  F \cup uS\not\in\Ti_{\GI,u}}
z^{|F|}\frac{F_{\VU'\setminus(S\cup V_F)}(z)}{{F}_{\VU'}(z)}\,.
\end{aligned}
$$

\\Therefore the function  $R^u_\GI(z)$ can be written as follows
\be
\label{Rugiz3}
\begin{aligned}
R^u_\GI(z) & = \sum_{S \subseteq \G_\GI(u) \atop |S|=1} z^{|S|}
+
\sum_{S \subseteq \G_\GI(u) \atop |S|=  2,~{\rm indep.}} z^{|S|}  -\sum_{S \subseteq \G_\GI(u) \atop |S|=  2,~{\rm indep.}} z^{|S|}\sum_{F \in \Fi_{\VU',S} \atop  F \cup uS\not\in\Ti_\GI}
z^{|F|}\frac{F_{\VU'\setminus(S\cup V_F)}(z)}{{F}_{\VU'}(z)}.\\
\end{aligned}
\ee
\\Now, by   Lemma \ref{lem37}, we have
\be\label{modK}
\left|z^{|F|}\frac{F_{\VU'\setminus(S\cup V_F)}(z)}{{F}_{\VU'}(z)}\right|= |z|^{|F|} \left| \frac{F_{\VU'\setminus(S\cup V_F)}(z)}{F_{\VU'}(z)} \right|\le |z|^{|F|} (1-a)^{-|V_F \cup S|}\le
\left( \frac{|z|}{1-a}\right)^{|F|}  (1-a)^{-|S|},
\ee
\vv
\\where the last inequality in \reff{modK} follows since
$|V_F \cup S| \leq |F| + |S|$ for all $F \in \mathcal{F}_{\VU',S} $. Hence, from \reff{Rugiz3} and \reff{modK} we get
\be
\begin{aligned}
|R^u_\GI(z)|
& \leq \sum_{S \subseteq \G_\GI(u) \atop |S|=1} |z|
+
\sum_{S \subseteq \G_\GI(u) \atop |S|=  2,~{\rm indep.}}|z|^{2} + \sum_{S \subseteq \G_\GI(u) \atop |S|=  2,~{\rm indep.}}|z|^{2}\sum_{F \in \Fi_{\VU',S} \atop  F \cup uS\not\in\Ti_\GI}
\left( \frac{|z|}{1-a}\right)^{|F|}  (1-a)^{-2}\\
& = \sum_{S \subseteq \G_\GI(u) \atop |S|=1} |z|
+
|z|^{2}
\sum_{S \subseteq \G_\GI(u) \atop |S|=  2,~{\rm indep.}}
\left(1 + \left({1\over 1-a}\right)^{2}\sum_{F \in \Fi_{\VU',S} \atop  F \cup uS\not\in\Ti_\GI}
\left( \frac{|z|}{1-a}\right)^{|F|} \right) \\
&\le \sum_{S \subseteq \G_\GI(u) \atop |S|=1} |z|
+
|z|^{2}
\sum_{S \subseteq \G_\GI(u) \atop |S|=  2,~{\rm indep.}}
\left(1 + \left({1\over 1-a}\right)^{2}\max_{S\in \G_\GI(u)\atop |S|=2,\;{\rm indep.}} \sum_{F \in \Fi_{\VU',S} \atop  F \cup uS\not\in\Ti_\GI}
\left( \frac{|z|}{1-a}\right)^{|F|} \right)\\
&=\D|z|+|I_u||z|^2\left(1 + \left({1-a}\right)^{-2}\max_{S\in \G_\GI(u)\atop |S|=2,\;{\rm indep.}}\sum_{F \in \Fi_{\VU',S} \atop  F \cup uS\not\in\Ti_\GI}
\left( \frac{|z|}{1-a}\right)^{|F|} \right)\\
& \le \D |z|+\k{\D^2\over 4}|z|^{2}\left( 1+ (1-a)^{-2}\max_{S\in \G_\GI(u)\atop |S|=2,\;{\rm indep.}} \sum_{F \in \Fi_{\VU',S} \atop  F \cup uS\not\in\Ti_\GI}
\left( \frac{|z|}{1-a}\right)^{|F|}  \right).\\
\end{aligned}
\ee
where in the last line we have used, via Definition \ref{double} and formula \reff{kaa}, that
$$
\sum_{S \subseteq \G_\GI(u) \atop |S|=2,~{\rm indep.}}1=|I_u|\le \k\left\lfloor{\D^2\over 4}\right\rfloor \le\k{\D^2\over 4}
$$
Now, setting  $|z|=(1-a){x\over \D}$, we get the following upper bound on $|R^u_\GI(z)|$ as $x\in [0,{1\over 2}]$

\be\label{eq:RP1P2}
\begin{aligned}
|R^u_\GI(z)| & \le (1-a)x+\k{x^2\over 4}\Big[(1-a)^2+\max_{S\in \G_\GI(u)\atop |S|=2,\;{\rm indep.}}
\sum_{F \in \Fi_{\VU',S} \atop  F \cup uS \not\in\Ti_\GI} \left({x\over \D}\right)^{|F|}\Big].
\end{aligned}
\ee
Let us now  bound
\be\label{fact}
\sum_{F \in \Fi_{\VU',S} \atop  F \cup uS \not\in\Ti_\GI} \left({x\over \D}\right)^{|F|}
\ee
\\Now fix an independent set $S=\{v_1,v_2\}\subseteq \G_\GI(u)$ {with, say, $v_1\prec v_2$}.  A forest $F\in\Fi_{\VU',S}$ that contributes to the sum above must satisfy $F \cup uS \notin \Ti_\GI$ and this can occur in two situations:
\begin{description}
    \item [Case 1: $F$ is a single tree.]
    In this case $F=\{\t\}$ where $\t$ is a Penrose tree containing both $v_1$ and $v_2$.
    Since $v_1$ and $v_2$ are independent, at least one additional edge is needed to connect them inside $\t$, hence $|F|\ge 2$.

    \item [Case 2: $F$ consists of two trees.]
    Here $F=\{\t_1,\t_2\}$ where each $\t_j$ is a Penrose tree containing (and rooted at) $v_j$, $j=1,2$. If $\GI \in \GG_0$ (claw-free but not necessarily square-free), by item c) of Lemma \ref{lem:case2-obstruction},
    { the tree $\t_1$} must be nonempty. If $\GI \in \GG_1$ (claw-free and square-free), by Corollary \ref{cor:one-trivial} both $\t_1$ and $\t_2$ must be nonempty.
\end{description}

\\Therefore, any such forest $F$ necessarily contributes at least one nontrivial tree rooted at $v_1$ and one tree (possibly trivial when $\GI\in \GG_0$) rooted at $v_2$, and its total weight is bounded by the product of the corresponding generating functions. We can thus bound \reff{fact} as
\begin{equation}\label{eq:fff}
 \sum_{\substack{F \in \Fi_{\VU',\{v_1,v_2\}} \\ F \cup \{u,v_1\}\cup \{u,v_2\} \notin \Ti_\GI}}
  {\left({x\over \D}\right)}^{|F|}
   \;\le\;
   \Bigg(\sum_{\substack{\t\in \Ti_{\GI',v_1} \\ |\t|\ge 1}} \left({x\over \D}\right)^{|\t|}\Bigg)
   \Bigg(\sum_{\substack{\t\in \Ti_{\GI',v_2} \\ |\t|\ge i}}\left({x\over \D}\right)^{|\t|}\Bigg).
\end{equation}
where $i=0,1$.

\\Recalling \reff{tgv}, for each $v_j=1,2$, we have
\[
T_{\GI',v_j}\left(y\right) \;=\; \sum_{\t\in \Ti_{\GI',v_j}} y^{|\t|}
= 1 + \sum_{\substack{\t\in \Ti_{\GI',v_j}\\ |\t|\ge 1}} y^{|\t|},
\]
where the term $1$ corresponds to the empty tree.
Therefore, for $i=0,1$,
since the degrees of $v_1$ and $v_2$ in $\GI'$ are at most $\D-1$, by Lemma \ref{le2} and recalling \reff{ginf2}, we obtain
\[
\sum_{\substack{\t\in \Ti_{\GI',v_j}\\ |\t|\ge i}}
\left({x\over \D}\right)^{|\t|}
\;\le\;\;g_{\D}\left({x\over \D}\right)-i\le h(x)-i
\]
and we can bound \reff{fact} as
\begin{equation}
 \sum_{F \in \Fi_{\VU',S} \atop  F \cup uS \not\in\Ti_\GI} \left({x\over \D}\right)^{|F|}
   \;\le\;
   \left(h(x)-1\right)
   \left(h(x)-i\right).
\end{equation}
for any independent pair $S=\{v_1, v_2\}\subset\G_\GI(u)$ and any $x\in[0,{1\over 2}]$.

\\In conclusion, plugging this estimate into \eqref{eq:RP1P2}, we find that for any
$z\in \C$ such that
$|z|=(1-a){x\over \D}$ with $x\in [0,{1\over 2}]$,
\be\label{rux}
\begin{aligned}
|R^u_\GI(z)| &\le (1-a)x+{\k}{x^2\over 4}\Big[(1-a)^2+ \left(h(x)-1\right)
   \left(h(x)-i\right)\Big]\\
 &= K^i_{\k}(a,x)
\end{aligned}
\ee
Therefore, by choosing $x = x^i_\k(a)$, we obtain that whenever
\[
   |z| \;\leq\; {(1-a)x^i_\k(a)\over \D},
\]
it holds that $|R^u_\GI(z)| \leq a$.
Since in this region we also have $F_{\VU-u}(z)\neq 0$ by the induction hypothesis, it follows that
\[
   F_\VU(z) \;\neq\; 0.
\]
Thus, the bound on $|R^u_\GI(z)|$ ensures the nonvanishing of $F_\GI(z)$ within the prescribed region, which
completes the proof of Proposition~\ref{pro1}.~~~~~~~~~~~~~~~~~~~~~~~~~~~~~~~~~~~~~~~~~~~~~~~~~~~~~~~~~~~~~~~~~~~~~~~~~~~~~~$\Box$

\subsubsection*{Acknowledgements}
P.M.S.F. was supported by
FAPEMIG (Funda\c{c}\~ao de Amparo \`a Pesquisa do Estado de Minas Gerais).
\subsubsection*{Competing Interests Statement}
The authors declare none.
\subsubsection*{Data Availability Statement}
No data was used for the research described in the article.
\subsubsection*{Funding statement}
This work was supported by  Funda\c{c}\~ao de Amparo \`a Pesquisa do Estado de Minas Gerais (FAPEMIG)
under research grants  30105/APQ-04481-22. The funder had no role in study design, data collection and analysis, decision to publish, or preparation of the manuscript.

\renewcommand{\theequation}{B.\arabic{equation}}
\setcounter{equation}{0}

\end{document}